\documentclass[12pt]{article}

\usepackage{latexsym}
\usepackage[centertags]{amsmath}
\usepackage{amsfonts}
\usepackage{amssymb}
\usepackage{amsthm}
\newcommand{\f}{\displaystyle \frac}

\begin{document}

\begin{center}

{ \Large \bf On the convergence of continued
fractions at Runckel's  points and the Ramanujan
conjecture}

\vspace{0.5cm}

Alexei Tsygvintsev\\

\vspace{0.5cm}

8 December 2004

\vspace{0.5cm}

U.M.P.A, Ecole Normale Sup\'erieure de Lyon\\
46, all\'ee d'Italie, F--69364  Lyon Cedex 07 \\
E-mail: atsygvin@umpa.ens-lyon.fr

\end{center}

\vspace{0.5cm}

\vspace{0.5cm}

\noindent{\bf  Abstract:} We consider  the limit
periodic continued
fractions of Stieltjes
$$
\frac{1}{1-}\, \, \frac{g_1 z}{1-}\, \, \, \,
\frac{g_2(1-g_1) z}{1-}\,\,
\frac{g_3(1-g_2)z}{1-\cdots,}, \, z\in \mathbb C,\,
g_i\in(0,1)\,, \lim\limits_{i\to \infty} g_i=1/2, 
\quad (1)
$$
appearing as Shur--Wall $g$-fraction representations
of certain
analytic self maps  of the unit disc $|w|< 1$, $w \in
\mathbb C$.
We precise the convergence behavior and prove the
general convergence  [2, p. 564 ] of (1) at the
Runckel's points of the singular line $(1,+\infty )$
It is  shown that in some cases  the  convergence
holds in the
classical sense. As a result a counterexample  to the
Ramanujan
conjecture [1, p.  38-39]  stating the divergence of a
certain
class of limit periodic continued fractions is
constructed.

\vspace{0.5cm}

{\Large \bf 1. Introduction}

\vspace{0.5cm}

 Let  $\mathcal E$ be  the class of  analytic self
maps of the unit disc
$\mathbb D=\{ w\in \mathbb C \,:\, |w|<1\}$. For 
$e\in \mathcal
E$ we introduce, following Shur [5],  the parameters
$\{ t_i
\}^{\infty}_{i=0}$ of $e(w)$ as follows
$$
e_0(w)=e(w),\quad t_0=e_0(0), \quad
e_{n+1}(w)=\displaystyle
\frac{1}{w}\displaystyle
\frac{e_n(w)-t_n}{1-\overline{t}_ne_n(w)},\quad
t_{n+1}=e_{n+1}(0)\,.
$$
The recursively defined functions
$$
[w;t_i]=t_i,\quad [w;t_l,\dots,t_k]=\displaystyle
\frac{t_l+w[w;t_{l+1},\dots,t_k]}{1+\overline{t}_lw[w;t_{l+1},\dots,t_k]},\quad
0\leq l<k\,, 
$$
provide then representation of $e(w)$
$$e(w)=\lim\limits
_{n\to \infty} [w;t_0,\dots,t_n]\,,$$ which
convergence  uniformly
over every compact  subset of $\mathbb D$.

 Let
$$
s_i(w;t)=t_i+\displaystyle
\frac{(1-|t_i|^2)w}{\overline{t}_iw+t^{-1}}\,, $$
 be an infinite
sequence of M\"{o}bius transformations of the variable
$t$. We
define $S_p(w;t)=s_0 \circ s_1\circ \cdots \circ
s_p(w;t)$. The
following formula can be easy derived:
$[w;t_0,\dots,t_n]=S_{n-1}(w;t_n)$, $n\geq 1$ and thus
$e(w)=\lim\limits _{n\to \infty} S_{n-1}(w;t_n)$.

Conversely with each sequence of complex numbers $t_i$
with
$|t_i|<1$, $i\geq 0$ one can associate an analytic
function $e(w)$
in $\mathbb D$ with $\sup\limits_{|w|<1}|e(w)|\leq 1$ 
such that
$t_i$ are just  as above.

We denote by $\mathcal W$ the set of continued
$g$-fractions of
Stieltjes
$$
g(z)=\frac{1}{1-}\, \, \frac{g_1 z}{1-}\, \,
\frac{g_2(1-g_1)z}{1-}\,\,
\frac{g_3(1-g_2)z}{1-\cdots}, \, z\in \mathbb C,\, 
g_i\in (0,1)\,. \eqno (1.1)
$$

\vspace{0.5cm}

 \noindent {\bf Definition 1.1.} {\it  Let
$h_i(z;t)=\displaystyle \frac{a_i}{1+t}$, $i\geq 0$ be
an infinite
sequence of M\"{o}bius transformations of the variable
$t$ with
$a_0=1$, $a_1=-zg_1$, $a_i=-g_i(1-g_{i-1})z$, $i\geq
2$ and let
$H_n(z;t)=h_0\circ h_1\circ \cdots \circ h_n(z;t)$.
Then the
continued fraction (1.1) is  called convergent at the
point  $z\in
\mathbb C$ if the limits
$$
\lim\limits _{n\to \infty}H_n(z;0)=\lim\limits _{n \to
\infty}
H_n(z;\infty)\,,
$$
exist in the extended complex plane $\mathbb
C_{\infty}$.}

\vspace{0.5cm}

As shown in [6, p. 279] the  continued fraction (1.1)
converges
for all $z$ from the  domain $C=\mathbb C_-\cup
\mathbb
C_+\cup(-\infty,1)$ to an analytic function $g(z)$
with  the property 
$\mathrm{Re}(\sqrt{1+z}\,g(z))>0$, $\forall \,z\in C$.

The Definition 1.1  of convergence is a classical one
and is rather unsatisfactory as one can imagine
situations then
$H_n(z;t)$ may converge at many $t$, but perhaps not
when $t$ is
$0$ or $\infty$ (see  [2, p. 565]). The following
refined
definition of convergence is due to Lisa Jacobsen [3,
p. 480]

\vspace{0.5cm}

\noindent {\bf Definition 1.2.} {\it The continued
fraction (1.1)
converges generally to a value $\alpha\in \mathbb
C_{\infty}$ at
$z \in \mathbb C_{\infty}$  if there exist sequences
$u_n$ and $v_n$ in $
\mathbb C_{\infty}$ such that
$$
\lim\limits _{n\to \infty}H_n(z;u_n)=\lim\limits _{n
\to \infty}
H_n(z;v_n)=\alpha, \quad \lim \inf_{n\to
\infty}\sigma(u_n,v_n)>0\,,$$ where $\sigma(x,y)$ is a
chordal
distance between $x,y\in \mathbb C_{\infty}$.}

\vspace{0.5cm}

One can  see  arbitrary
M\"{o}bius transformations $l_n(z)$ as isometries of 
the
hyperbolic space $\mathbb H^3$ [2, p. 559]. Then the
above
definition  is fully justified by the
following geometric result due to  Beardon

\vspace{0.5cm}

\noindent {\bf Theorem 1.3.}  ([2, p. 567]) {\it A
sequence $l_n$
of M\"{o}bius maps converges generally to $\alpha\in
\mathbb
C_{\infty}$ iff $l_n\to \alpha$ pointwise on $\mathbb
H^3$.}

\vspace{0.5cm}

A continued fraction convergent  in the
classical sense always  converges generally to the
same value  (one puts  $u_n=0$, $v_n=\infty$, $n\geq 0
$ in Definition 1.2). The
next theorem describes the  correspondence between
functions from
the classes $\mathcal E$ and $\mathcal W$

\vspace{0.5cm}

\noindent {\bf Theorem 1.4.} ([6, p. 289]) {\it Let
$$
w=-1+2(1-\sqrt{1-z})/z, \eqno (1.2)
$$ be the conformal mapping of $C$ onto
$\mathbb D$ with a positive square root branch for 
$z<1$. To
every function $e(w)\in \mathcal E$ corresponds  a
function
$g(z)\in \mathcal W$ according to
$$
\displaystyle \frac{1+w}{1-w}\displaystyle \frac{1-w
e(w)}{1+we(w)}=g(z)\,,
$$
where the coefficients $t_i$ and $g_i$ are related by
$t_{k-1}=1-2g_k$, $k\geq 1$.}

\vspace{0.5cm}

{\Large \bf 2. The Runckel's points}

\vspace{0.5cm}

Let $t_i$, $i=0,1,\dots$ be a sequence of real 
numbers with $|t_i| < 1$. We
assume that 
$$
\sum_{i=0}^{\infty}t_{i}^2<\infty\,, \eqno(2.1)
$$
and therefore 
$$
\lim\limits_{i\to \infty} t_i=0 \,. \eqno (2.2)
$$

Let $E \subset \mathcal E$ be the subset of functions
$e(w)$
whose parameters $t_i$ satisfy the above conditions.

\vspace{0.5cm}

\noindent {\bf Definition 2.1.}  {\it  A point $r\neq
\pm 1$, $|r|=1$
  is a Runckel's  point for $e(w)\in E$ if
the limit $e(r)=\lim\limits _{n\to \infty}
[r;t_0,\dots,t_n]$
exists and is equal to $1$.}

\vspace{0.5cm}

We note that, according to Runckel  [4, p. 98],  if in
addition to (2.1), (2.2)  there exists a natural  $p$
such that $\sum\limits_{i=0}^{\infty}
|t_{i+p}-t_i|<+\infty$, then, as $k\to \infty$,
$[w;t_0,\dots,t_k]$ converges uniformly over every
compact subset
of $|w|\leq 1$, $w^p\neq 1$ to $e(w)$ (analytic in
$\mathbb D$)   continuous  and $|e(w)|<1$ for all $w$
in $|w|\leq 1$, $w^p\neq
1$. Thus, in this particular case, every  Runckel's 
point $r$ satisfies    $r^p=1$. We will see some
examples of these functions in Section 3.

The following  result concerns the general convergence
 of 
$g(z)$  at Runckel's points

\vspace{0.5cm}

\noindent {\bf Theorem 2.2.} {\it Let $e(w)\in E$, let
$r$ be
its Runckel's point  and $g(z)\in \mathcal W$ be the
corresponding
$g$-fraction given by Theorem 1.4. Then $g(z)$
converges generally
at the point
$$
z_r=2(1+\mathrm{Re}(r))^{-1}>1\,, \eqno (2.3)
$$ which is the image of $r$ by the conformal mapping
(1.2) (called  also a Runckel's point
of $g(z)$).}

 \vspace{0.5cm}

\noindent {\bf Proof.} The following formula relates
the partial
approximants of the fractions $e(z)$ and $g(z)$  (see
Theorem
78.1, [6])

$$
\begin{array}{lll}
\displaystyle \frac{1+w}{1-w}\displaystyle
\frac{1-wS_n(w;t)}{1+wS_n(w;t)}=H_{n+1}(p(w);l_{n+1}(w;t)),
\\ \\

H_{n+1}(p(w);l_{p+1}(w;t))=h_0\circ h_1(p(w);t)\circ
h_2(p(w);t)\circ\cdots  \\ \\ \cdots \circ
h_{n+1}(p(w);l_{n+1}(w;t)), \quad  n \geq 1\,,
\end{array} \eqno (2.4)
$$
where
$l_{n+1}(w;t)=-2\displaystyle\frac{(1-g_{n+1})(1-t)w}{(1-wt)(1+w)}$,
$p(w)=4w^2/(1+w)^2$.

Let $r$ and $z_r$  be the Runckel's points of $e(w)$
and $g(z)$
respectively as given in  (2.3). We note that
 the parameters  $t_{k-1}=1-2g_k$ are real thus 
$e(w)$ always has two
complex Runckel's points $r$ and $\overline{r}$,
$r\neq
\overline{r}$ corresponding to unique
$z_r=p(r)=p(\overline{r})>1$
 belonging to the singular line of the continued
fraction $g(z)$.
Putting $t=t_{n+1}$, $w=r,\overline{r}$, $z=z_r$ in
(2.4) and
taking the limit as $n \to \infty$ we obtain 
$$
\lim\limits _{n\to
\infty}H_{n+1}(z_r,u_{n+1})=\lim\limits _{n \to
\infty} H_{n+1}(z_r,v_{n+1})=1\,, \eqno (2.5)
$$
with $u_{n+1}=\overline{ v_{n+1}
}=l_{n+1}(r;t_{n+1})$. Finally
one observes that $\lim \inf\limits _{n\to
\infty}\sigma(u_n,v_n)>0$ that leads to the general
convergence of
$g(z)$ at $z=z_r$.  $\Box$

The sequence $H_n(z_r;0)$, $n\geq 1$ can be
still divergent in $\mathbb C_{\infty}$ i.e $g(z)$
being divergent
at $z=z_r$ in the classical sense. Nevertheless we
have the
following result 

\vspace{0.5cm}

\noindent { \bf Lemma 2.2. } {\it One has
$$\mathrm{min}\{|H_n(z_r;0)-1|,|H_{n+1}(z_r;0)-1|\}\to
0, \quad \mathrm{as} \quad    n\to
\infty\,. \eqno (2.6)
$$}

\noindent {\bf Proof.}

To see it we follow the  Beardon's geometric method
described in
[2]. Let $\gamma_{0,\infty}$ be the vertical geodesic
in $\mathbb
H^3$ with endpoints $0$ and $\infty$, let
$\gamma_{u_n,v_n}$ be
the geodesic with endpoints $u_n$ and $v_n$. One
checks that the
hyperbolic distance $d$ between $\gamma_{0,\infty}$
and
$\gamma_{u_n,v_n}$ always satisfies $d>0$. Since
$H_n(z_r;t)$ are
isometries of $\mathbb H^3$, the distance between
$H_n(z_r;\gamma_{0,\infty})$ and
$H_n(z_r;\gamma_{u_n,v_n})$ is
also $d$. According to (2.5), as $n\to \infty$ the
geodesic
$H_n(z_r;\gamma_{u_n,v_n})$ shrinks to the point $1\in
\mathbb C_{\infty}$ and thus
$\mathrm{min}\{|H_{n+1}(z_r;0)-1|,|H_{n+1}(z_r;\infty)-1|\}\to
0$ as $n\to \infty$
that  gives (2.6).  $\Box$

We denote  $\mathrm{dist}(X,Y)$ the Euclidean distance
between two subsets $X,Y\subset \mathbb C_{\infty}$.

The following theorem describes the possible limit
behavior of the sequence $H_n(z_r;0)$ as $n\to
\infty$.

\vspace{0.5cm}

\noindent {\bf Theorem  2.3.}  $\lim \limits_{n\to
\infty} \mathrm{dist}(H_n(z_r;0),\{1,\displaystyle
\left(\frac{1+r}{1-r}\right)^2\})=0$

\noindent {\bf Proof.} Putting  $t=1$, $w=r$  in (2.4)
one obtains
$$
\displaystyle \frac{1+r}{1-r}\displaystyle
\frac{1-r[r;t_0,t_1,\dots,t_n,1]}{1+r[r;t_0,t_1,\dots,t_n,1]}=H_{n+1}(z_r;0)\,.
\eqno (2.7)
$$
We define $A_{\nu}$, $C_{\nu}$, $\tilde A_{n}$,
$\tilde C_{n}$ by the following formulas
$$
\begin{array}{ll}
A_0=1,\quad C_0=t_0, \\ \\
A_{\nu+1}=A_{\nu}+t_{\nu+1}r^{\nu+1}\bar{C}_{\nu}, \\
\\ C_{\nu+1}=C_{\nu}+t_{\nu+1}r^{\nu+1}\bar{A}_{\nu},
\quad 0\leq \nu \leq n-1,\\ \\
\tilde A_{n}=A_{n}+r^{n+1}\bar{C}_{n},  \quad \tilde
C_{n}=C_{n}+r^{n+1}\bar{A}_{n}\,.
\end{array} \eqno (2.8)
$$

Then, as follows  from    [4, p. 99],  we have 
$$
[r;t_0,t_1,\dots,t_{\nu}]=\f{ C_{\nu} } { A_{\nu} },
1\leq \nu \leq n, \quad [r;t_0,t_1,\dots,t_{n},1]=\f{
\tilde C_{n} } { \tilde A_{n}}\,,
$$
and  
$$
1-\left |\f{ C_{n} }{ A_{n} } \right |^2=\f{ P_{n}} {
| A_{n} |^2 }\,. \eqno (2.9)
$$

$$
\f{ \tilde C_{n}}{\tilde A_{n}}-\f{C_{n}}{A_{n}}=\f{
P_{n} }{A_{\nu}\tilde A_{n}}\,, \eqno (2.10)
$$
where $P_{n}=\prod\limits_{k=0}^{n}(1-t_k^2)$.

According to  Shur  [6, p. 138]    
$$
\tilde C_{n}=r^{\nu+1}\bar{\tilde A}_{n}\,, \eqno
(2.11)
$$
that gives in particular $|\tilde C_{n} |=| \tilde
A_{n}|$, $n\geq 1$.

In view of  $\lim\limits_{n \to \infty}C_{n} /
A_{n}=1$  the formula (2.9) implies that
$\lim\limits_{n \to \infty} A_{n}=\lim\limits_{n \to
\infty} C_{n}=\infty$)   and $\lim\limits_{n \to
\infty}(A_n-C_n)=0$ since $P_n$ 
is bounded from zero by (2.1).
We have also  $ \tilde A_{n},\tilde C_{n} \neq 0$,
$\forall \,n$ by $|C_n|<|A_{n}|$. 

If  $\tilde A_{n}$ is bounded from zero as $n \to
\infty$, then  $\lim\limits_{n \to \infty}\tilde C_{n}
/ \tilde A_{n}= 1$ according to  (2.10)  and   the
Theorem 2.3   is proved. Let now assume that there
exists an infinite  sequence of $\nu$,
$\nu_1<\nu_2<\cdots$ such that
$\lim \limits_{i\to \infty} \tilde A_{\nu_i }=0$.
Then,  as seen from (2.8) and $\lim\limits_{n \to
\infty}C_{n} /  A_{n}=1$ ,   as $i \to \infty$,
$\arg(A_{\nu_i})$ is close to $(\pi+\arg(r^{\nu_i+1}))
/ 2$, $ \mod \,( \pi)$ and $\arg(\tilde A_{\nu_i})$ is
close to 
one of the following values: $\{\arg(A_{\nu_i})\pm
\pi/2,\arg(A_{\nu_i}),\arg(A_{\nu_i})+\pi\}$.  By
(2.11) that implies $\lim \limits_{i\to \infty}
\mathrm{dist}(\tilde C_{\nu_i} / \tilde
A_{\nu_i},\{1,-1\})=0$.  The result follows then from
(2.7). $\Box$

\vspace{0.5cm}

{\Large \bf 3. The  Ramanujan's conjecture on the
convergence of
limit periodic continued fractions}

\vspace{0.5cm}

It is interesting to find  examples of functions
$e(w)\in E$ with the  classical convergence at the
Runckel's points.

We define 
$e_p(w)=(1+w^p) / 2$, $p\in \mathbb N$. Then, as
shown in [5, p. 142] (see also [4, p. 106]), 
$$
\begin{array}{ll}
e_p(w)=[w;t_0,0,\dots,0,t_1,\dots,0,t_2,0,\dots], \\ 
\\ t_0=1/2,
\quad t_n=2/(2n+1), \quad n\geq 1\,,
\end{array} 
$$
where $p-1$ zeros are added between $t_n$ and
$t_{n+1}$.  According to Shur the function 
 $e_p(w)\in E$  is continuous in
$|w|\leq 1$ with Runckel's points $w$ given by the
roots of  $w^p=1$, $w\neq \pm 1$. The
coefficients $g_i$ of the corresponding to  $e_p(w)$
continued
fraction $g_p(z)$ are defined by relations
$g_k=(1-t_{k-1})/2$ and 
$$
g_1=\frac{1}{4},\quad g_{1+l}=\frac{1}{2}, \quad l\not
\equiv
0\,\, \mathrm{mod}\,  (p), \quad
g_{1+pk}=\displaystyle
\frac{2k-1}{2(2k+1)}, \quad k\geq 1\,.
$$
Introducing the  parameters $b_1=g_1$,
$b_i=g_i(1-g_{i-1})$, $i\geq
2$ one obtains
$$
\begin{array}{ll}
b_1=\displaystyle\frac{1}{4}, \quad
b_2=\displaystyle\frac{3}{8},
\quad b_{1+l}=\displaystyle\frac{1}{4}, \quad l\not
\equiv 0,1 \,
\, \mathrm{mod}\, (p),\quad
b_{1+pk}=\displaystyle\frac{1}{4}\frac{2k-1}{2k+1},\\
\\
b_{2+pk}=\displaystyle \frac{1}{4}\frac{2k+3}{2k+1},
\quad k\geq
1\,.
\end{array}
$$
In particular  $\lim\limits _{i\to \infty} b_i=1/4$.

The continued fractions  $g_p(z)$ has  the Runckel's
 point  
$z_{p}= \cos^{-2}( \pi / p)>1 $
corresponding to $r=\exp(2\pi i /p)$ -- the Runckel's
point of $e_p(w)$.

Repeating the arguments used in the proof of Theorem
2.3 it is easy to show that   for every  odd $p>1$ the
limit
periodic continued fractions
$$
g_p(z_{p})=\frac{1}{1-}\, \, \frac{b_1z_{p}}{1-}\, \,
\frac{b_2z_{p}}{1-}\,\,\frac{b_3z_{p}}{1-\cdots},
\quad \lim \limits_{i\to
\infty}b_iz_{p}=z_{p}/4>1/4\,,
$$
 converges  to $1$ and is divergent by oscillations 
if $p$ is even.
 
  Indeed, as seen from formulas (2.8), the
coefficients  $A_{\nu}$, $C_{\nu}$ are always real and
positive for $\nu \geq 0$ and thus  
 $\tilde A_{\nu}=A_{\nu}+r^{\nu+1}\bar{C}_{\nu}$,
$r=\exp(2\pi i /p)$  is bounded or not from $0$ as
$\nu \to \infty$ depending on whether $p$ is odd or
even.

  In the  simplest case $p=3$ the above
 continued fraction takes the form
$$
\frac{1}{1-}\,\,\frac{1}{1-}\,\,\frac{3/2}{1-}\,\,\frac{1}{1-}\,\,\frac{1/3}{1-}\,\,\frac{5/3}{1-}\,\,\frac{1}{1-}\,\,
\frac{3/5}{1-}\,\,\frac{7/5}{1-}\,\,\frac{
1}{1-}\,\,\frac{5/7}{1-}\,\,\frac{9/7}{1-}\,\,\frac{1}{1-}\,\,
\frac{7/9}{1-}\,\,\frac{11/9}{1-}\,\,\frac{1}{1-\cdots}\,,
$$
and converges to $1$.

For a general limit periodic continued fraction
$$
\frac{1}{1-}\, \, \frac{a_1}{1-}\, \,
\frac{a_2}{1-\cdots,}\,,
$$
given by a sequence of real parameters $a_i$, $i\geq
1$ such that
$\lim\limits_{i\to \infty} a_i =a>1/4$ Ramanujan
stated that
it is always divergent (see [1, p. 38-39] for
details).  The above continued fraction represents a 
counterexample to this conjecture.

\vspace{0.5cm}

{\Large \bf 3. Acknowledgments.}

\vspace{0.5cm}

The author is grateful to A. Glutsyuk for useful
discussions.

\vspace{0.5cm}

{\Large \bf List of references}

\vspace{0.5cm}

\noindent [1] G. E. Andrews,  B. C. Berndt, L.
Jacobson,  R.L. Lamphere, The  continued fractions
found in the  unorganized portions of Ramanujan's
notebooks,  
Mem. Amer. Math. Soc. 99, no. 477, 1992

\noindent [2]  A.F. Beardon, Continued fractions,
discrete groups
and complex dynamics, Computational Methods and
Function Theory,
Volume 1 (2001), No. 2, 535-594 

\noindent [3] L. Jacobsen, General convergence of
continued
fractions, Trans. Amer. Math. Soc. 294 (1986), 477-485

\noindent [4] H.-J., Runckel , Bounded analytic
functions in the
unit disk and the behavior of certain analytic
continued
fractions near the singular line,  J. reine angew.
Math. 281
(1976), 97-125

\noindent [5] I. Shur, \"{U}ber Potenzreihen, die im
Innern des
Einheitskreises beschr\"{a}nkt sind, J. reine angew.
Math. 147
(1916), 205-232, and 148 (1917), 122-145

\noindent [6] H. S. Wall, Analytic theory of continued
fractions,
D. van Nostrand Company Inc., NY, London, 1948

\end{document}